\documentclass{amsart}
\usepackage{latexsym}
\usepackage{amsmath}
\usepackage{amssymb}
\usepackage[all]{xy}
\usepackage{comment}

\usepackage{autobreak}
\usepackage{enumitem}
\usepackage[hypertexnames=false]{hyperref} 
\usepackage[capitalize]{cleveref}
\crefname{equation}{}{}
\crefname{enumi}{}{}

\usepackage[\ifdefined\enabletodonotes\else disable\fi]{todonotes}

\usepackage{autonum} 

\usepackage{tikz}
\usetikzlibrary{positioning}
\usetikzlibrary{arrows.meta}

\numberwithin{equation}{section}
\newtheorem{thm}{Theorem}[section]
\newtheorem{prop}[thm]{Proposition}
\newtheorem{cor}[thm]{Corollary}
\newtheorem{lem}[thm]{Lemma}

\theoremstyle{definition}
\newtheorem{defn}[thm]{Definition}

\newtheorem{assertion}[thm]{Assertion}
\theoremstyle{remark}
\newtheorem{rem}[thm]{Remark}
\newtheorem{ex}[thm]{Example}

\newcommand{\K}{{\mathbb K}}
\newcommand{\V}{{\mathcal V}}
\newcommand{\W}{\mathcal{W}}
\newcommand{\C}{{\mathcal C}}
\newcommand{\e}{\varepsilon}
\newcommand{\R}{{\mathbb R}}
\newcommand{\mapright}[1]{%
 \smash{\mathop{%
  \hbox to 1cm{\rightarrowfill}}\limits_{#1}}}
\newcommand{\maprightd}[2]{%
 \smash{\mathop{%
  \hbox to 0.5cm{\rightarrowfill}}\limits^{#1}\limits_{#2}}}
\newcommand{\mapleft}[1]{%
 \smash{\mathop{%
  \hbox to 1cm{\leftarrowfill}}\limits_{#1}}}
\newcommand{\mapleftu}[1]{%
 \smash{\mathop{%
  \hbox to 0.8cm{\leftarrowfill}}\limits^{#1}}}
\newcommand{\maprightu}[1]{%
 \smash{\mathop{%
  \hbox to 1cm{\rightarrowfill}}\limits^{#1}}}
\newcommand{\maprightud}[2]{%
 \smash{\mathop{%
  \hbox to 1cm{\rightarrowfill}}\limits^{#1}_{#2}}}
\newcommand{\mapleftud}[2]{%
 \smash{\mathop{%
  \hbox to 1cm{\leftarrowfill}}\limits^{#1}_{#2}}}

\newcommand{\modelcat}{\mathcal{M}}
\newcommand{\et}[1]{\mathsf{et}({#1}^{[0, \infty)})}

\author[K. Kuribayashi]{Katsuhiko Kuribayashi}
\address{%
  Department of Mathematical Sciences,
  Faculty of Science,
  Shinshu University,
  Matsumoto, Nagano 390-8621, Japan
}
\email{kurimath@shinshu-u.ac.jp}

\begin{document}
\title[Enriched nerves, interleaving distances and the multiplicative MPSS]{Enriched nerves, interleaving distances and the multiplicative magnitude-path spectral sequence 
}


\begin{abstract} By regarding the nerves of  categories enriched over a monoidal poset as persistence simplicial sets, we investigate an upper bound of  the interleaving distance in the homotopy category between the nerves in the sense of Lanary and Scoccola \cite{LaSc}. In particular, it holds that  the distance between the nerves for $r$-homotopy equivalent generalized metric spaces is less than or equal to the half of $r$. 
Moreover, a product on the cochain complex associated with a persistence simplicial set is introduced and related to the cup product on the magnitude cochain complex in the sense of Hepworth \cite{H}. As a consequence, the product gives the multiplicative structure to the cohomological magnitude-path spectral sequence (MPSS) whose first term is isomorphic to the magnitude cohomology ring. Furthermore,  we discuss at which term the MPSS  for a connected finite digraph collapses by using 
the interleaving of the enriched nerve of the digraph. 
 \end{abstract}

\maketitle


\section{Introduction}
{\it Magnitude homology} for graphs is introduced by Hepworth and Willerton \cite{HW}. Leinster and Shulman \cite{LS_21} genelarize the magnitude homology to a large class of categories enriched over a monoidal category. Subsequently, Otter \cite{O} has given a framework to relate magnitude homology with persistence homology which appears in topological data analysis. 
The key idea is to introduce a simplicial set indexed by the poset $[0, \infty]$ for a given $[0, \infty]^{\rm op}$-category $\mathcal{C}$, namely a generalized metric space. The persistence simplicial set is called the {\it enriched nerve} of $\mathcal{C}$. 
The {\it blurred magnitude homology} due to Otter, which is a persistence version of magnitude homology, is further studied in 
\cite{C, GH, Me}. Very recently, Ivanov and Xu 
 \cite{IX} have given a lot of fascinating results on $\nu$-{\it Vietoris-Rips simplicial sets} by which both the blurred magnitude chain complex and the Vietoris-Rips complex are induced. 
 To compare such persistence objects homotopically, we may utilize the notion of the interleavings up to homotopy which is originally introduced for persistence spaces by Blumberg and Lesnick \cite{Bl-L} and is investigated in a more general setting by Lanari and Scoccola \cite{LaSc}. 
 
It is worthwhile mentioning that there exists a spectral sequence, that is  called the magnitude-path spectral sequence \cite{HW, Asao, DOMZ} (MPSS for short), converging to the {\it reachability homology} of a directed graph (digraph) $G$ with the shortest path metric; see \cite{HR}. The first page of the MPSS is indeed the magnitude homology of $G$ and the second one contains the {\it path homology} of $G$ in the sense of 
Grigor'yan, Lin,  Muranov and Yau \cite{GLMY} on the horizontal line.  This intimate relationship of the two homologies has been clarified by Asao \cite{Asao}. 
Moreover, the entire first page is the {\it bigraded path homology}  due to Hepworth \cite{HR_2026} of the given digraph. 

The aim of this article is to consider the enriched nerves in the category of persistence simplicial sets endowed with the projective model structure. Furthermore, we discuss a product on the cochain complex associated with the enriched nerve which gives a multiplicative structure to the cohomological MPSS.
 
More precisely, by regarding enriched nerves as persistence objects, we give an upper bound of the {\it interleaving distance $d_{\rm IHC}$ in the homotopy category} between such nerves introduced by Lanari and Scoccola \cite{LS}.  In particular, Theorem \ref{thm:interleavings} asserts that the distance $d_{\rm IHC}$ between the enriched nerves for $2r$-homotopy equivalent generalized metric spaces is less than or equal to $r$.  
As a corollary, we see that the half of the {\it diameter} of a generalized metric space $M$ is 
an upper bound of the distance $d_{\rm IHC}$ between the enriched nerves for $M$ and a point; see Corollary \ref{cor:diameter}.  
As seen in Remark \ref{rem:homotopy},  the colimit of the enriched nerve of every connected digraph is homotopy equivalent to that of a point. However, the distance $d_{\rm IHC}$ indeed detects {\it evolution} of the nerves by the persistence parameter; 
see Corollaries \ref{cor:toy} and \ref{cor:toy2} as toy computational examples of $d_{\rm IHC}$. 

In other words, the distance $d_{\rm IHC}$ gives a lower bound for the $r$-homotopy equivalence of two given generalized metric spaces $X$ and $Y$. Thus, for example, if the {\it bottleneck distance} of the persistent homology of the enriched nerves for $X$ and $Y$ is greater than or equal to $s$, then there is no $2r$-homotopy equivalence between $X$ and $Y$ for $r < s$; see Remark \ref{rem:BD}. 
While the interleaving distance of the homology groups of the enriched nerves of metric spaces is considered in \cite[Section 6]{Me} and \cite[Section 4]{IX}, 
our objects in the study are persistence simplicial sets themselves.  

In the meddle of  this article, 
we introduce a product on the dual to the chain complex associated with a persistence simplicial set, 
which is compatible with persistence structure on the dual; see Theorem \ref{thm:main_I}.  In fact, the product is obtained by restricting the cup product on the dual to the reachability chain complex discussed in \cite{HR}. 
As a consequence, the product induces that on the magnitude cochain complex in the sense of Hepworth \cite{H}; see Proposition \ref{prop:MC}. We observe that the magnitude cohomology algebra is a complete isometry invariant for extended quasi-metric spaces $(X, d)$ each of which $\inf \{ d(a, b) \mid a, b \in X, a\neq b\}$ is positive; see \cite[Theorem 3.1 (Recovery theorem)]{H}. 
Moreover, the product in Theorem  \ref{thm:main_I} allows us to give a multiplicative structure to the cohomological MPSS converging to the reachability cohomology algebra of a given metric space; see Theorem \ref{thm:SS}.   

The MPSS is constructed by a persistence simplicial chain complex whose domain is the poset of non-negative integers.   
Such persistence objects are incorporated into the notion of {\it extended tame} persistence objects; see \cite{CGL, KNSWY}. 
In the latter half of this article, 
we consider at which term the MPSS for a connected finite digraph collapses by using the interleaving distance $d_{\rm IHC}^{\rm \ \! et}$ in the homotopy category of extended tame chain complexes between the  enriched nerve complexes of the digraph and a point; see Theorem \ref{thm:IHC_MPSS}. Thus, the distance $d_{\rm IHC}^{\rm \ \! et}$ from the acyclic chain complex to the enriched nerve complex of  a digraph $G$ gives a sufficient condition for $G$ to be $r$-{\it contractible} in one of the Brown categories of $\mathbb{N}$-metric spaces in the sense of Livernet, Roff and  Whitehouse; see \cite[Section 6]{LRW} and also \cite[Section 7]{HR_2026}.

The rest of this manuscript is organized as follows. 
In Section \ref{sect:Nerves}, following Otter \cite{O}, we introduce the enriched nerve for a $[0, \infty]^{\rm op}$-category. 
Section \ref{sect:E_nerves} deals with interleaving distances of enriched nerves and proves Theorem \ref{thm:interleavings}. In Section \ref{sect:MPSS}, we address 
the multiplicative MPSS and prove Theorem \ref{thm:SS}.  In Section \ref{sect:ID_MPSS}, after recalling the model structure of extended tame persistence objects which takes values in a model category, Theorem \ref{thm:IHC_MPSS} is proved. In Appendix \ref{sect:et_proj}, we compare the distances  $d_{\rm IHC}$ and $d_{\rm IHC}^{\rm \ \! et}$ for persistence objects whose domains are  the poset of non-negative integers; see Proposition \ref{prop:et_proj} and Corollary \ref{cor:N(X)}.

\section{The nerve of the category enriched over a monoidal category}\label{sect:Nerves} 
In what follows, we may call a functor from a poset  $\V$ to a category a {\it persistence object}, for example, a {\it persistence simplicial set} for  a functor $\V \to \mathsf{sSet}$. 

We recall the nerve of an enriched category introduced by Leinster, Shulman \cite{LS} and Otter \cite{O}. Let $\V=(\V, \otimes, 1_\V)$ be a monoidal category and $C$ a small 
$\V$-category; that is,  
$C$ has a set of objects and an object $C(x, y)$ in $\V$ is assigned for elements $x$ and $y$ in $C$ together with a composition $C(x,y)\otimes C(y ,z) \to C(x, y)$ and identity $1_{\V} \to C(x, x)$ which satisfy associativity and unitalily conditions\footnote{We refere the reader to \cite[PART I, Section 3]{Riehl} for basic concepts of enriched categories.}. 
By definition, a morphism $f : C \to C'$ of $\V$-categories is given by an object map $x \mapsto f(x) \in C'$ for $x \in C$ together with morphisms $f_{x,y} : C(x, y) \to C'(f(x), f(y))$ in $\V$ for $x, y \in C$ such that the diagrams  
\[
\xymatrix@C15pt@R20pt{
1_{\V}  \ar[r] \ar[rd] & C(x, x) \ar[d]^{f_{x,x}} &  \! \! \! \! \! {\rm and}  \! \! \! \! \!  & C(y, z)\otimes C(x, y) \ar[r]^-{\circ} \ar[d]_{f_{y,z}\otimes f_{x,y}}  & C(x, z) \ar[d]^{f_{x,z}} \\
& C'(f(x), f(x)) &   &  C'(f(y), f(z))\otimes C'(f(x),  f(y)) \ar[r]_-{\circ}  & C'(f(x), f(z)) 
}
\]
commute for $x, y, z \in C$. We have the category $\V\text{-}\mathsf{Cat}$ of $\V$-categories.

Assume further that $\V$ is {\it semicartesian}; that is, the unit $1_\V$ is the terminal object. Then, a functor 
$N(C) : \V^{\rm op} \to  \mathsf{sSet}$
which is called the {\it enriched nerve} of $C$ is defined by 
\[
N(C)(\e)_n := \sqcup_{x_0, x_1, ..., x_n \in \rm{Ob}C}\V(\e, C(x_0, x_1)\otimes C(x_1, x_2)\otimes \cdots \otimes C(x_{n-1}, x_n)), 
\]
where the face maps $d_i :  N(C)(\e)_n \to N(C)(\e)_{n-1}$ and the degeneracy maps  $s_j :  N(C)(\e)_n \to N(C)(\e)_{n+1}$ are defined using composition morphisms 
and identity-assigning morphisms, respectively. Observe that $d_0$ and $d_n$ are defined by using the unique map 
$C(x, y) \to 1_\V$. 

In the case where $(\V, \leq, \otimes, 1_\V)$ is a monoidal poset whose opposite category $\V^{\rm op}$ is semicartesian and $C$ is a $\V^{\rm op}$-category, we see that 
\[
N(C)(\e)_n = \{(x_0, x_1, ..., x_n) \mid x_i \in {\rm Ob}C , \ell(x_0, x_1, ..., x_n) \leq \e \  {\rm in \ } \V\}, 
\]
where  $\ell(x_0, x_1, ..., x_n)$ denotes the element  $C(x_0, x_1)\otimes C(x_1, x_2)\otimes \cdots \otimes C(x_{n-1}, x_n)$ in $\V$. We call the element $\ell(x_0, x_1, ..., x_n)$ the {\it length} of the tuple $(x_0, ..., x_n)$. 
Thus, the nerve construction  above gives a functor $N : \V^{\rm op}\text{-}\mathsf{Cat} \to  \mathsf{sSet}^{\V}$ from categories of $\V^{\rm op}$-categories to the category of  persistence simplicial sets. 
We may write $cN(C)$ for the simplicial set $\text{colim}_lN(C)(l)_\bullet$ obtained by the colimit construction. 

\begin{rem}\label{rem:ceiling-F}
Let $F : \V^{\rm op} \to \mathcal{M}^{\rm op}$ be a lax monoidal functor between monoidal category. For a $\V^{\rm op}$-category $C$, one has a 
$\mathcal{M}^{\rm op}$-category $FC$ whose composition $FC(a, b)\otimes FC(b, c) \to FC(a, c)$ is defined by the composite 
\[
\xymatrix@C45pt@R15pt{
Fx\otimes Fy \ar[r]^-{\alpha_{x, y}} & F(x\otimes y) \ar[r]^-{F(x\otimes y \to z)} & Fz
}
\]
with the structure map $\alpha_{x,y}$ of $F$ and the composition $x\otimes y \to z$, where $x = C(a, b)$, $y= C(b, c)$ and $z = C(a, c)$. 
For example, the ceiling function $\lceil \ \rceil : [0, \infty] \to \overline{\mathbb{N}}:=\mathbb{N}\cup \{\infty\}$ gives rise to a lax monoidal functor  $[0, \infty)^{\rm op} \to \mathbb{N}^{\rm op}$
between the monoidal posets. Observe that $\lceil a \rceil + \lceil b \rceil \geq \lceil a +b  \rceil$. 
\end{rem}

\begin{ex}\label{ex:L} i) (An enriched nerve of a metric space)
Let $(X, d)$ be a generalized metric space (an extended quasi-psudometric space)\footnote{A Riemannian diffeological spaces \cite{KSS} containing Riemannian manifolds and a small full subcategoery of persistence objects with the interleaving distance \cite{KNSWY} are such spaces. In particular, the homotopy set of topological spaces is a generalized metric space; see \cite[Theorem 5.2]{KNSWY}.}; that is, the distance function $d : X\times X \to [0, \infty]$ satisfies the conditions: 
\begin{enumerate}
\item[i)] $d(x, x) = 0$ and 
\item[ii)] $d(x, y) \leq d(x, z) + d(z, y)$. 
\end{enumerate}
Then, the space $X$ is identified with a $[0, \infty]^{\rm op}$-category defined by $X(x, y) = d(x. y)$ for $x, y \in X$; see  \cite{Law}. 
Applying the enriched nerve construction mentioned above,  we have a persistence simplicial set $N(X) : [0, \infty] \to \mathsf{sSet}$ of $X$ with the form 
\[
N(X)(l)_n : = \big\{(x_0, x_1, ... x_n) \mid x_i \in X, \sum_{i=0}^{n-1} d(x_i, x_{i+1}) \leq l \big\}
\]
whose face maps $d_i $ and degeneracy maps $s_j$ are defined by  
\[
d_i ((x_0, x_1, ..., x_n))= (x_0, ..., x_{i-1}, x_{i+1}, ..., x_n)
\]
for $i = 0, ..., n$ and 
\[
s_j ((x_0, x_1, ..., x_n))= (x_0, ..., x_j, x_j, ..., x_n)
\]
for $j = 0, ..., n$, respectively.  By definition, the length $\ell(x_0, ..., x_n)$ has the form  
$\ell((x_0, ..., x_n)) =  \sum_{i=0}^{n-1} d(x_i, x_{i+1})$. We see that the construction gives rise to a functor 
$N : [0, \infty]^{\rm op}\text{-}\mathsf{Cat} \to \mathsf{sSet}^{([0, \infty], \leq)}.$ 

Moreover, the nerve yields a persistence chain complex $\K N(X)$ whose homology is isomorphic to the blurred magnitude homology in the sense of Otter; see \cite[Definition 9.3 and Theorem 9.4]{O}. 
The colimit $\text{colim}_l \K N(X)(l)$ is nothing but the {\it reachability complex} $RN(X)$ which induces the magnitude-path spectral sequence; see Section \ref{sect:MPSS_1}.  

ii) Let $X$ be a compact Riemannian manifold, more general, a compact Riemannian diffeological space defined in \cite{KSS}. We have an 
$\mathbb{N}^{\rm op}$-category of the type 
$\lceil \ \rceil X$ and the a persistence simplicial set $N(\lceil \ \rceil X)$; see Remark \ref{rem:ceiling-F}. 

\end{ex}

 \section{Interleavings between enriched nerves}\label{sect:E_nerves} Throughout this section, 
 let $\V$ be the monoidal poset $([0, \infty),  \leq,  +, 0)$.  
 Observe that $\V^{\rm op}$ is semicartesian. Thus,  we have the functor $N : \V^{\rm op}\text{-}\mathsf{Cat} \to  \mathsf{sSet}^{\V}$. 
 For a generalized metric space $X$, the enriched nerve $N(X)$ defined in Example \ref{ex:L} is a persistence simplicial set. 
The functor $N : [0, \infty]^{\rm op}\text{-}\mathsf{Cat} \to \mathsf{sSet}^{([0, \infty], \leq)}$ described in Example \ref{ex:L} i) and the inclusion functor $\iota: \V \to [0, \infty]$ give the functor $\iota^*N : [0, \infty]^{\rm op}\text{-}\mathsf{Cat} \to \mathsf{sSet}^\V$ defined by $(\iota^*N)(C)(l)= N(C)(\iota(l))$. 
 The functor will be also denoted by $N$. In what follows, we denote by $\W$  the monoidal poset $\V$ or $[0, \infty]$. 
 
 Some of results in this section are generalized to the case of a more general monoidal poset. While it is crucial to address for example the monoidal poset $(\mathbb{R}^n, \leq, +, 0)$ with the product order $\leq$ in the study of multiparameter  persistence theory, we do not pursue such generalizations in this manuscript except for the argument in Remark \ref{rem:a_generalization}. 
 
 \subsection{Interleavings and $r$-homotopies}\label{sect:I-r}
 In this section, we consider the interleaving distance between persistence simplicial sets $N(C)$ 
 associated with $\W^{\rm op}$-categories via the functor $N :  \W^{\rm op}\text{-}\mathsf{Cat} \to  \mathsf{sSet}^{\V}$.   
 We begin by reviewing the interleaving distance introduced in \cite{CCGGO} and results in \cite{B-S} related to the distance, which we use extensively in this article.
By regarding the poset $(\V, \leq)$ as a category, we have the functor category $\C^{(\V, \leq)}$ for a category $\C$.
For a real number $\e \geq 0$, define a functor $T_\e : (\V, \leq) \to (\V, \leq)$ by $T_\e(a) = a+ \e$.  Moreover, the $\e$-{\it shift functor}
\[
( \ )^\e :  \C^{(\V, \leq)} \to \C^{(\V, \leq)}
\]
is defined by $( \ )^\e (F) = F^\e := FT_\e$.
We have a natural transformation $\eta^\e : id_{\C^{(\V, \leq)}} \Rightarrow (\ )^\e$ defined by the structure maps 
$\eta^\e_{F}(i):= F(i \leq i + \e)  : F(i) \to F(i+\e)=F^\e(i)$.

\begin{defn}\label{defn:interleaving}
(\cite[Definition 4.2]{CCGGO}, \cite[Definition 3.1]{B-S})
Objects $F$ and $G$ in $\C^{(\V, \leq)}$ are {\it $\e$-interleaved} if there exist natural transformations $\varphi : F \to GT_\e$ and
$\psi : G \to FT_\e$
such that the diagram
\begin{eqnarray}\label{eq:interleaving_1}
\xymatrix@C35pt@R20pt{
F \ar[r]^{\eta^\e_F}  \ar[dr]^(0.7){\varphi} &  F^{\e}  \ar[r]^{\eta^\e_{F^\e}} \ar[rd]^(0.7){\varphi^\e}|\hole & F ^{\e+\e}   \\
G  \ar[r]_{\eta^e_G}  \ar[ur]^(0.3){\psi}|\hole& G^\e \ar[ur]^(0.3){\psi^\e} \ar[r]_{\eta^\e_{G^\e}}  & G^{\e+\e}
}
\end{eqnarray}
is commutative. 
Such a tuple $(F, G, \varphi, \psi)$ is called an {\it $\e$-interleaving}.
\end{defn}

For objects $F$ and $G$ in $\C^{(\V, \leq)}$, 
the  {\it interleaving distance} between $F$ and $G$ is define by
\[
d_{\text{I}}(F, G) \! :=\! \text{inf}(\{ \e \geq 0 \mid \text{$F$, \!$G$ are $\e$-interleaved} \}\cup \{\infty\}).
\]
We observe that $F$ and $G$ are $\e+\delta$-interleaved for each $\delta \geq 0$ if the objects are $\e$-interleaved; see 
\cite[Lemma 3.4]{B-S}.

 In what follows, we use terminology in \cite{DS, H_Model} for model categories. 
 Let $\mathcal{M}$ be a cofibrantly generated model category and
$\mathcal{M}^{(\V, \leq)}$ the model category endowed with the {\it projective model structure} whose weak equivalences and fibrations are pointwise weak equivalences and pointwise fibrations of $\mathcal{M}$, respectively; see \cite[Theorem 11.6.1]{H_Model}. 

We denote by $\mathsf{Top}$ and $\mathsf{sSet}$ the cofibrantly generated model categories of topological spaces and simplicial sets described for example in \cite[11.1.5 Example]{H_Model}, respectively.


We say that objects $F$ and $G$ in $\mathcal{M}^{(\V, \leq)}$ are $\e$-{\it interleaved in the homotopy category} if they are
$\e$-interleaved in $\text{Ho}(\mathcal{M}^{(\V, \leq)})$. The shift functor $( \ )^\e : \mathcal{M}^{(\V, \leq)} \to \mathcal{M}^{(\V, \leq)}$ preserves weak equivalences. Thus,
we can consider the commutative diagram (\ref{eq:interleaving_1}) for $F$ and $G$  in $\text{Ho}(\mathcal{M}^{(\V, \leq)})$. 
For objects $F$ and $G$ in $\text{Ho}(\mathcal{M}^{(\R, \leq)})$, Lanari and Scoccola \cite{LaSc} introduce the {\it interleaving distance in the homotopy category}\label{index:homotopyID} define by
\[
d_{\text{IHC}}(F, G) \! :=\! \text{inf}(\{ \e \geq 0 \mid \text{$F$, \!$G$ are $\e$-interleaved in the homotopy category} \}\cup \{\infty\}).
\]

By the same definition as above, we define the {\it interleaving distance in the homotopy category} for the objects $F$ and $G$ in 
$\text{Ho}(\mathcal{M}^{(\V, \leq)})$. 

\begin{rem}\label{rem:IHC}
An interleaving in $\text{Ho}(\mathcal{M}^{(\V, \leq)})$ mentioned above is considered with a homotopy commutative diagram of the form (\ref{eq:interleaving_1})
in $\mathcal{M}^{(\V, \leq)}$ using fibrant-cofibrant replacements $RQF$ and $RQG$ of $F$ and $G$, respectively.  In fact, by definition, the Hom-set $\text{Hom}_{\text{Ho}(\mathcal{M}^{(\V, \leq)})}(F, G)$ between $F$ and $G$ is given by the set 
of homotopy classes of morphisms from $RQF$ to $RQG$: 
\[
{\rm Hom}_{\text{Ho}(\mathcal{M}^{(\V, \leq)})}(F, G) = {\rm Hom}_{(\mathcal{M}^{(\V, \leq)})}(RQF, RQG)/ \simeq.
\] 
We require that the homotopy which gives the homotopy commutative diagram is a morphism in $\mathcal{M}^{(\V, \leq)}$. 
\end{rem}
 
In what follows, let $\K$ be a commutative ring with unit. 
 For a  functor $P : \V\to \mathsf{sSet}$ to the category of simplicial sets, we have a functor 
$\K P :  \mathcal{V} \to \mathsf{Ch}_\K$ to the category of chain complexes 
defined by $(\K P(l), d_n = \sum_{i=0}^n (-1)^{i}d_i)$ for $l \geq 0$,  
where $\K P(l)_n$ is the free $\K$-module generated by the set $P(l)_n$  
and $d_i$ are the face maps in $P(l)$. 
Let $H_n$ denote the homology functor $\mathsf{sSet}^\V \to \mathsf{Ab}^\V$ defined by $H_n(P) := H_n(\mathbb{Z}P, d_*)$, 
where  $\mathsf{Ab}$ is the category of abelian groups. It follows from \cite[Proposition 3.6]{B-S} that 
\begin{eqnarray}\label{eq:I-IHC}
d_{\rm I} (H_n(X_\bullet), H_n(Y_\bullet)) \leq d_{\rm IHC} (X_\bullet, Y_\bullet)
 \end{eqnarray}
 for persistence simplicial sets $X_\bullet$ and  $Y_\bullet$ and for each $n\geq 0$. 
 
 In the rest of this section, we consider an upper bound of the interleaving distance in the homotopy category between enriched nerves.

 \begin{defn}\label{defn:r-homotopy}(cf. \cite[\S 2.2]{I}, \cite[Definition 5.2]{IX}, \cite[Definition 2.2.2]{LRW}, \cite[Definition 2.2]{HR}) Let $r$ be an object in $\V$ and $f, g : C \to C'$ morphisms of $\W^{\rm op}$-categories, namely $1$-Lipschitz maps in $\W$; see the beginning of Section \ref{sect:Nerves}. We say there is  an {\it elementary $r$-homotopy} from $f$ to $g$, and write  $f  \rightsquigarrow_r g$, if $C'(f(x), g(x))) \leq r$ in $\V$ for every $x \in C$. 
 
Moreover, maps $f$ and $g$ are $r$-{\it homotopic}, denoted $f\sim_r g$,  if there is a zig-zag sequence of elementary $r$-homotopies; that is, there exists a tuple $f_1 , ..., f_k$ of  morphisms of $\V^{\rm op}$-categories such that $f= f_1$ and $f_k = g$ and, each $i$, either $f_i  \rightsquigarrow_r f_{i+1}$ or  $f_{i+1}  \rightsquigarrow_r f_i$. By definition, a morphism $f : C\to C'$ is an $r$-homotopy equivalence if there exists a morphism $g :C' \to C$ such that $g\circ f \sim_r id_C$ and $f\circ g \sim_r id_{C'}$. 
 \end{defn}
 
 \begin{rem}\label{rem:homotopy} (1) Let $f : C\to C'$ be a morphism of $\W^{\rm op}$-categories.  Then, the nerve functor $N$ mentioned in Section \ref{sect:Nerves} gives a morphism $N(f)_r : N(C) \to N(C')^r$ defined by the composite $\eta^r(N(C')) \circ N(f) : N(C) \to N(C') \to N(C')^r$
 for an object $r$ in $\V^{\rm op}$; see the sentence before Definition \ref{defn:interleaving} for $\eta^r$.

(2)  Morphisms $f, g : C \to C'$ with $f  \rightsquigarrow_r g$ give rise to the simplicial homotopy 
 $\{ \{h_j : N(C)(l)_n \to N(C')(l+r)_{n+1} \}_{0\leq j\leq n}\}$ from $N(f)_r$ to $N(g)_r$\footnote{$d_0h_0=N(f)_r$ and $d_{n+1}h_n=N(g)_r$. } defined by 
 \[
 h_j(x_0, ..., x_n)= (g(x_0), ..., g(x_j), f(x_j), ..., f(x_n)).
 \] 
 As a consequence, we see that $cN(C):=\text{colim}_lN(C)(l)_\bullet$ is homotopic to $cN(pt)$ if there exists $r_C$ in $\V$ such that $C(x, y) \leq r_C$ in $\V$ for $x$ and $y$ in $C$.

(3)  Let $f, g : K \to L$ be simplicial maps and $\{\{h_j : K_n \to L_{n+1} \}_{0\leq j\leq n}\}$ a simplicial homotopy from $f$ to $g$. Then, we have a simplicial map 
 $H : K\times \Delta[1] \to L$ with $H\circ \e_0 = f$ and $H\circ \e_1 = g$. In fact, the homotopy $H$ is defined by 
\begin{eqnarray}\label{eq:H}
H_n(x_n, c_j) = 
\left\{
 \begin{array}{ll}
 g_n(x_n) & \text{for $j= 0$}, \\
  f_n(x_n) & \text{for $j = n+1$},  \\ 
  d_j h_{j-1} &\text{for $1 \leq j\leq n$}.
  \end{array}
  \right. 
 \end{eqnarray}
 Here, we regard $\Delta[1]_n = \Delta([n], [1])$ as the set  $\{ c_0, ..., c_{n+1}\}$ with $c_i^{-1}(0) = \{0, ..., i-1\}$; see, for example, 
 \cite[Proposition 10.5.2]{R}. 
  \end{rem}
 
 \begin{prop}\label{prop:r-homotopy} {\rm (cf. \cite[Proposition 3.3]{IX})} With the same notation as in Remark  \ref{rem:homotopy},  
 for morphisms $f, g : C \to C'$  of $\W^{\rm op}$-categories with $f  \rightsquigarrow_r g$, the maps $|N(f)_r|$ and $|N(g)_r|$ are homotopic 
 in the model category $\mathsf{Top}^{(\V, \leq)}$ endowed with the projective model structure.   
  \end{prop}

To prove Proposition \ref{prop:r-homotopy}, we need the following two lemmas whose proofs are postponed to  the end of this section. 

\begin{lem}\label{lem:cof}
The object $N(C)(*)_\bullet^\e$ is cofibrant in $\mathsf{sSet}^\V$ for each $\e\geq 0$.
\end{lem}

\begin{lem}\label{lem:key1}
Let $X$ be a cofibrant object in $\mathsf{Top}^{(\V, \leq)}$. Suppose that $H: X\times I \to Y$ is a morphism with $\e_0H = f$ and $\e_1H =g$, where 
$I$ is the constant persistence space defined by $I(l) =[0, 1]$ for $\l \in \V$ and $(\e_i)(l)(x) = (x , i)$ for $x \in X(l)$.  Then, $f$ is right homotopic to $g$ and hence 
$f\simeq g$ in $\mathsf{Top}^{(\V, \leq)}$.
\end{lem}

\begin{proof}[Proof of Proposition \ref{prop:r-homotopy}] 
By \cite[Theorem 11.6.5]{H_Model} and Lemma \ref{lem:cof}, we see that the persistence space $|N(C)(*)_\bullet^\e|$ is cofibrant. 
Observe that $|N(C)(*)_\bullet^\e| = |N(C)(*)_\bullet|^\e$ and $|N(C)(*)_\bullet|^\e$ is fibrant; see Remark \ref{rem:IHC}. 

The maps $N(f)$ and $N(g)$ are regarded as 
maps $cN(C) \to cN(C')$. 
Then, Remark \ref{rem:homotopy} allows us to obtain a commutative diagram 
\begin{eqnarray}\label{eq:diagram_IHC}
\xymatrix@C30pt@R15pt{
cN(C) \times \Delta[1] \ar[r]^H & cN(C') \\
N(C)(l)_\bullet \times \Delta[1] \ar[r]_{H_{{\rm res} ,l}} \ar[u] \ar[d]  & N(C')(l+ r)_\bullet  \ar[u] \ar[d]\\
N(C)(l')_\bullet \times \Delta[1] \ar[r]_{H_{{\rm res}, l'}}  \ar@/^3.5pc/[uu]& N(C')(l'+ r)_\bullet  \ar@/_3.5pc/[uu],
}
\end{eqnarray}
where $H_{{\rm res}, k}$ is the restriction to $N(C)(k)_\bullet$ of $H$ and 
the vertical arrows denote the inclusions induced by the structure maps of the persistence simplicial set $N(C)(\ast)_\bullet$.
Thus, we see that the realizations $|H_{\rm res}|=\{ |H_{{\rm res}, l}| \}_{l \in \V}$ of the homotopies $H_{{\rm res}, l}$ gives rise to a homotopy from 
$|N(f)_r|$ to $|N(g)_r|$ in $\mathsf{Top}^{(\V, \leq)}$. In fact, the commutativity of the lower square in (\ref{eq:diagram_IHC}) enables us to deduce that 
the map $|H_{\rm res}| : |N(C)(*)_\bullet| \times I \to |N(C)(*)^r_\bullet|$ is in $\mathsf{Top}^{(\V, \leq)}$.  
Lemma \ref{lem:key1} gives the homotopy we require. 
\end{proof} 

\begin{rem}
A cofibration in  $\mathcal{M}^{(\V, \leq)}$ endowed with the projective model structure is complicated. Indeed, it is a retract of a transfinite composition; see \cite[Theorem 11.6.1]{H_Model} for more details. Thus, in order to prove Proposition \ref{prop:r-homotopy},  we use the notion of right homotopy 
instead of that of left homotopy.  
\end{rem}

\begin{thm}\label{thm:interleavings}
Let $f : C \to C'$ be a $2r$-homotopy equivalence in $\W^{\rm op}$-categories. Then, the persistence simplicial sets 
$N(C)$ and $N(C')$ in $\mathsf{sSet}^{(\V, \leq)}$ are $r$-interleaved in the homotopy category. As a consequence, one has 
\[
d_{\rm IHC}(N(C), N(C')) \leq r.
\]
\end{thm}

\begin{proof} There exists a morphism $g : C'\to C$ such that $f\circ g \sim_{2r} id$ and $g\circ f \sim_{2r} id$. 
By applying Proposition \ref{prop:r-homotopy} repeatedly, we see that $|(N(g)_r)^r|\circ |N(f)_r| \simeq |\eta_{N(C)}^{r+r}|$ and 
 $|(N(f)_r)^r|\circ |N(g)_r| \simeq |\eta_{N(C')}^{r+r}|$. 
 In fact, it holds that 
 \[
 (N(g)_r)^r\circ N(f)_r = N(g\circ f)_{r+r} \  \  \text{and}  \ \ \eta^{r+r}_{N(C)} = N(id)_{r+r}.
 \]
 This yields that $d_{\rm IHC}(S_*(|N(C)|), S_*(|N(C')|))\leq d_{\rm IHC}( |N(C)|,  |N(C')|) \leq r$, where $S_*( \ )$ denotes the singular simplex functor. Since $S_*(|F|)\cong F$ in ${\rm Ho}(\mathsf{sSet}^{(\V, \leq)})$ for a cofibrant persistence simplicial set $F$, we have the result; see \cite[Theorems 8.5.23 and 11.6.5]{H_Model}.
\end{proof}

\begin{rem}\label{rem:a_generalization}
Let $\V$ be a monoidal poset. The notion of an $r$-homotopy in Definition \ref{defn:r-homotopy} is generalized to that for objects in 
$\V^{\rm op}$-category. Moreover, the proof of Theorem  \ref{thm:interleavings} allows us to conclude that the assertion of the first half 
of the theorem holds replacing the term ``a $2r$-homotopy equivalence" with ``an $(r+r)$-homotopy equivalence" for 
$\mathcal{W}=\V$. 
\end{rem}

\begin{cor}\label{cor:diameter} Let $x_0$ be an object in $C$. Suppose that the supremum $r_C$ 
 of the set $\{ C(x, x_0) \mid x \in {\rm ob} \ C \}$\footnote{If $C$ is a generalized metric space, the supremum is less than or equal to the diameter of $C$.} exists. 
 Then, 
 the objects  
$N(C)_\bullet,$ and  $N(\{x_0\})_\bullet$ 
in $\mathsf{sSet}^{(\V, \leq)}$
are  $\frac{r_C}{2}$-interleaved in the homotopy category. 
 \end{cor}

\begin{proof}
The inclusion from the trivial category $\{x_0\}$ to $C$ is an $r_C$-homotopy equivalence. Then, the result follows from Theorem \ref{thm:interleavings}. 
\end{proof}


\begin{rem}\label{rem:BD}
For a finite generalized metric space $A$, the $k$th homology group $H_k(\mathbb{Z}N(A)(l) ; \K)$ with coefficients in a field $\K$ is of finite dimension for each $k$ and $l$.  Let $\mathcal{B}_{A}$ denote 
the barcode of the persistence module $H_k(\mathbb{Z}N(A)(*); \K)$. 
Then, by \cite[Proposition 3.6 and Theorem 4.16]{B-S} and the inequality (\ref{eq:I-IHC}), we see that, for a finite generalized metric spaces $X$ and $Y$,  the bottleneck distance 
$
d_{\rm B}(\mathcal{B}_{X}, \mathcal{B}_{Y})
$ 
of the persistent homology gives a lower bound of 
the interleaving distance in the homotopy category $d_{\rm IHC}(N(X), N(Y))$. 

If the bottleneck distance is greater than or equal to $s$, then 
$X$ and $Y$ is {\it not} $2r$-homotopy equivalent for any $r$ less than $s$.  This follows from Theorem \ref{thm:interleavings}. 
We observe that $H_k(\mathbb{Z}N(X)(*); \mathbb{Z})$ is isomorphic to the blurred magnitude homology; see Example \ref{ex:L} i).

We refer the reader to \cite[Theorem 3.3]{KNWY} for the equivalence of the distance $d_{\rm IHC}$ of chain complexes with other interleaving distances up to homotopy; see \cite{Bl-L, LaSc} for interleavings up to homotopy. 
\end{rem}

\subsection{Toy examples as corollaries of Theorem \ref{thm:interleavings}}
In this section, we give computational toy examples for $d_{\rm IHC}$.  

 \begin{cor}\label{cor:toy}
 Let $K_n$ be the complete graph endowed with the shortest path distance, where $n \geq 2$. 
 Then, one has 
 $\displaystyle{d_{\rm IHC}(N(K_n)_\bullet, N(pt)_\bullet) = \frac{1}{2}}$.
 \end{cor}
 
\begin{proof} Corollary \ref{cor:diameter} enables us to conclude that 
$d_{\rm IHC}(N(H_{m, n})_\bullet, N(x_1)_\bullet) \leq \frac{1}{2}$.
For each integer $k \geq 0$, 
the structure map $N(K_n)(k)_\bullet \to N(K_n)(k+\e)_\bullet$ is identity for $\e <1$. 
The same condition holds for $N(pt)_\bullet$.   
Therefore, if we have $d_{\rm IHC}(N(K_n)_\bullet, N(pt)_\bullet) < \frac{1}{2}$, then 
$d_{\rm IHC}(N(K_n)_\bullet, N(pt)_\bullet) =0$. Thus, there exists a weak  equivalence between $N(K_n)_\bullet$ and $N(pt)_\bullet$. 
In particular, the simplicial sets $N(K_n)(0)_\bullet= \coprod_n\Delta[0]$ and $N(pt)(0)_\bullet= \Delta[0]$ are weak homotopy equivalent. 
This is a contradiction.
\end{proof}


\begin{cor}\label{cor:toy2}
Let $(T_{m, n}, d)$ be the finite metric space consisting three points $x_0$, $x_1$ and $x_2$ with 
$d(x_0, x_1) = m$, $d(x_1, x_2) = n$ and $d(x_0, x_1) = m+n$, where $m\leq n$. Then, 
$\displaystyle d_{\rm IHC}(N(T_{m, n})_\bullet, N(x_1)_\bullet) = \frac{n}{2}$.
\end{cor}

\begin{proof}
We write $T$ for the metric space $T_{m,n}$.  
Corollary \ref{cor:diameter} enables us to deduce that 
$d_{\rm IHC}(N(T)_\bullet, N(x_1)_\bullet) \leq \frac{n}{2}$. 
For the case where $m=n$, the result follows from the same argument as in Corollary \ref{cor:toy}. 

Suppose  that $m < n$.  Assume further that $d_{\rm IHC}(N(T)_\bullet, N(x_1)_\bullet) < \frac{n}{2}$. It follows from the inequality (\ref{eq:I-IHC}) that 
\[
\displaystyle d_{\rm I}(H_0(N(T), H_0(N(x_1))) \leq d_{\rm IHC}(N(T)_\bullet, N(x_1)_\bullet) < \frac{1}{2}. 
\]
Thus, there exists a positive real number $\e$ such that $H_0(N(T))$ and $H_0(N(x_1))$ are $(\frac{n}{2}-\e)$-interleaved and $\frac{m}{2} < \frac{n}{2}-\e$.  
We have a commutative diagram 
\[
\xymatrix@C25pt@R15pt{
H_0(N(T)(0)) \ar[r]^{\theta} \ar[rd]_{\varphi} & H_0(N(T)(m)) \ar[r]^-{\cong}  & \cdot \ar[r] & H_0(N(T)(n)  \\
H_0(N(x_1)(0)) \ar[r]_-=  & H_0(N(x_1)(\frac{n}{2}-\e) \ar[ru]_{\psi} \ar[rr]_-= & &  H_0(N(x_1)(m)),
} 
\]
where $\theta$ denotes the structure map $H_0(N(H)(0 < m)$. 
We have $H_0(N(T)(0)) \cong \mathbb{Z}\{ x_0-x_1, x_1, x_2\}$. Since $\mathbb{Z}N(T)_0(m)$ is generated by 
$x_0, x_1, x_2$ and the component 
$\mathbb{Z}N(T)_1(m)$ is generated by  $(x_0, x_1)$ and $(x_1, x_0)$, it follows that 
$H_0(N(T)(m)) \cong \mathbb{Z}\{ x_1, x_2\}$ and the map $\theta$ is an epimorphism. 
The diagram above gives a commutative triangle 
\[
\xymatrix@C30pt@R10pt{
\mathbb{Z}^{\oplus 3} \ar@{->>}[rr] \ar[rd]_\varphi & & \mathbb{Z}^{\oplus 2}  \\
& \mathbb{Z} \ar[ru]_\psi
}
\]
in which the horizontal arrow is an epimorphism. Then, the map $\psi$ is an epimorphism, which is a contradiction. 
\end{proof}

\subsection{Proofs of Lemmas \ref{lem:cof} and \ref{lem:key1}}\label{sect:A1}

Before proving Lemma \ref{lem:cof}, we recall a result due to Lanari and Scoccola \cite{LaSc} which is  
a characterization of a cofibrant persistence simplicial set in 
$\mathsf{sSet}^{(L, \leq)}$ endowed with the projective model structure, where $(L, \leq)$ is an arbitrary poset. 

\begin{defn}\label{defn:functions} (\cite[Definition 4.1]{LaSc}) A {\it $L$-filtered simplicial set} $(X, \beta)$ is a simplicial set $X$ equipped with the  set of functions $\beta = \{\beta_n : X_n \to L\}_{n\geq 0}$, satisfying 
\begin{itemize}
\item $\beta_{n-1}(d_i(\sigma)) \leq \beta_n(\sigma)$ for every $n \geq 1$, $\sigma \in X_n$, and face map $d_i : X_n \to X_{n-1}$; 
\item $\beta_{n+1}(s_j(\sigma)) \leq \beta_n(\sigma)$ for every $n \geq 0$, $\sigma \in X_n$,  and degeneracy map $s_j : X_n \to X_{n+1}$. 
\end{itemize}
\end{defn}

\begin{defn} (\cite[Definition 4.2]{LaSc}) Given a $L$-filtered simplicial set $(X, \beta)$, define a persistence simplicial set 
$\widehat{(X, \beta)} \in \mathsf{sSet}^L$ such that for $r \in L$ we have 
\[
\widehat{(X, \beta)}(r)_n = \{\sigma \in X_n \mid \beta_n(\sigma)\leq r\}, 
\]
with faces and degeneracies given by restricting the ones of $X$.
\end{defn}

The following is the result of the characterization mentioned above. 

\begin{prop}\label{prop:cof} {\rm (\cite[Proposition 4.5]{LaSc})} A persistence simplicial set $X$ is isomorphic to one associated with a $L$-filtered simplicial set if and only if $X$ is cofibrant in $\mathsf{sSet}^L$.  
\end{prop}

\begin{proof}[Proof of Lemma \ref{lem:cof}]
It follows from Proposition \ref{prop:cof} that $N(C)(*)_\bullet$ is cofibrant in $\mathsf{sSet}^\V$. In fact, we have a $\V$-filtered  simplicial set $(cN(C), \beta)$, where $\beta_n(x_0 , ...,  x_n)= \ell(x_0 , ...,  x_n)$ with the length $\ell$; see the paragraph before Remark \ref{rem:ceiling-F} for the notations. Then, it follows that $\widehat{(cN(C), \beta)}= N(C)(*)_\bullet$. 

Define functions $\beta'_n : X_n:=cN(C)_n \to \V$ by $\beta'_n (\sigma) = \text{max}\{0, \beta_n(\sigma) -\e\}$. 
It is readily seen that the set of functions $\beta'=\{\beta'_n\}$ satisfy the two conditions in Definition \ref{defn:functions}. By Proposition \ref{prop:cof}, in order to  prove the result, it suffices to show that 
the third equality in 
\begin{align}
\widehat{(X, \beta)}^\e(r)_n = \widehat{(X, \beta)}(r+\e)_n &= \{\sigma \in X_n \mid \beta_n(\sigma) \leq r +\e \} \\
&=  \{\sigma \in X_n \mid \beta'_n(\sigma) \leq r \} = \widehat{(X, \beta')}(r)_n=:E_n
\end{align}
holds. The first, second and fourth equalities follow from the definitions. 

For $\sigma \in E_n$, suppose that $\beta'(\sigma) = 0$. Then 
$\beta_n(\sigma) -\e \leq 0$. Since $\beta_n(\sigma) \leq \e$, it follows that $\sigma \in  \widehat{(X, \beta)}(\e)_n \subseteq  \widehat{(X, \beta)}(r+\e)_n$. If $\beta'_n(\sigma) > 0$, then $\beta_n(\sigma) -\e = \beta'_n(\sigma) \leq r$. We see that $\sigma$ is in $\widehat{(X, \beta)}^\e(r)_n$. 
Let $\sigma$ be a simplex in $\widehat{(X, \beta)}^\e(r)_n$.  Suppose that  $\beta_n(\sigma) -\e \geq 0$. We have 
$\beta'_n(\sigma) = \beta_n(\sigma) -\e \leq r$. If $\beta_n(\sigma) -\e < 0$, then $\beta'(\sigma) = 0 \leq r$. This implies that 
$\sigma \in E_n$. 
\end{proof}

\begin{proof}[Proof of Lemma \ref{lem:key1}]
By considering the adjoint $ad(H)$ to $H$, we have a commutative diagram 
\begin{eqnarray}\label{eq:H}
\xymatrix@C30pt@R15pt{
X \ar[r]^{ad(H)} \ar[rd]_{(f, g)}& Y^I \ar@{->>}[d]^{\pi}\\
& Y\times Y,
}
\end{eqnarray}
where $\pi$ is pointwise evaluation fibration. Observe that $\pi$ is a fibration in $\mathsf{Top}^{(\V, \leq)}$. 
Let $c : Y\to Y^I$ be the map assigns the constant map at $y$ to a point $y \in Y$. 
Then, we have a commutative diagram 
\begin{eqnarray*}
\xymatrix@C40pt@R10pt{
Y  \ar@{>->}[rd]_\iota^\sim  \ar@/^1.2pc/[rrr]^{\Delta} \ar[rr]_-{c}^-\sim  & & Y^I \ar@{->>}[r]_-\pi & Y\times Y \\
& Z \ar@{->>}[ru]^p  & X \ar[u]_{ad(H)} \ar[l]^{K}&
}
\end{eqnarray*}
in which $p\circ i$ is a decomposition of $g$ and $K$ is a lift of $ad(H)$. We observe that $p$ is a weak equivalence and $X$ is  cofibrant. 

Since the diagonal map $\Delta$ factors through $Z$, it follows that $Z$ is a path object for $Y$. Then, the commutativity of the diagram (\ref{eq:H}) allows us to deduce that the morphism $K$ is a right homotopy between $f$ and $g$. 
Each object in $\mathsf{Top}^{(\V, \leq)}$ is fibrant. Then, it follows from \cite[Proposition 7.4.8]{H_Model}  that  
$f\simeq g$ in $\mathsf{Top}^{(\V, \leq)}$. 
\end{proof}

\section{A product on a cochain complex associated with an enriched nerve}\label{sect:MPSS}

 As in Section \ref{sect:E_nerves},  let $\V$ denote the monoidal poset $([0, \infty), \leq, +, 0)$. 
Let $P : \V\to \mathsf{sSet}$ \todo{Let $P : \mathcal{V} \to \mathsf{sSet}$ a functor from a modoidal poset.} be a persistence simplicial set and $\K P$ the persistence chain complex associated with $P$. Assume that 
the structure map $P(\e \leq \e')$ is injective for $\e \leq \e'$. We may identify the element $\sigma$ in $P(\e)$ with $P(\e\leq \e')\sigma$ in $P(\e')$. 

Let $cP$ and $P(<\e)$ denote the colimits $\text{colim}_l P(l)$ and $\text{colim}_{i< \e} P(l)$, respectively. We observe that 
$P(\e) = \text{colim}_{i\leq \e} P(l)$. It follows that the structure map $P(\e' < \e)$ induces a chain map $\iota_{< \e} : \K P(<\e) \to  \K P(\e)$. 
The chain complex  $RP:=\K cP\cong\text{colim}_l\K P(l)$ is called the {\it reachability (chain) complex}.  

\subsection{Magnitude homology and cohomology}\label{sect:MH} 
We recall magnitude cohomology introduced by Hepworth; see \cite[Section 2, Definition 5.3 and Example 5.6]{H}. 

Let $X$ be a generalized metric space. 
The (unnormalized) magnitude chain complex $MC_{\e, *}(X)$ of $X$ is defined by  
\[
MC_{\e, n}(X):=\K\{(x_0, ..., x_n) \mid \sum_{i=0}^n d(x_i, x_{i-1}) = \e\}
\]
with the differential $d_n :=\sum_{i=0}^n(-1)^id_i$ where for $\sigma = (x_0, ..., x_n)$, 
\[
d_i(\sigma)=\left\{
\begin{array}{ll}
 (x_0, ..., x_{i-1}, x_{i+1}, ..., x_n) & \text{if  $(x_0, ..., x_{i-1}, x_{i+1}, ..., x_n)$ is in  $MC_{\e, n-1}(X)$}  \\
 0 & \text{otherwise. }
\end{array}
\right.
\]
The {\it magnitude homology} $MH_{\e, *}(X)$ is the homology of the magnitude chain complex.  
The magnitude cochain complex of $X$ is the dual to the magnitude chain complex: 
$MC_{\e}^*(X) = \text{Hom}(MC_{\e, *}(X), \K)$. The {\it magnitude cohomology} $MH_{\e}^*(X)$ is the cohomology of the magnitude cochain complex. 

Give $\varphi \in MC_{\e}^m(X)$ and $\psi \in MC_{\delta}^n(X)$, the {\it product } $\varphi\cdot \psi \in MC_{\e+\delta}^{m+n}(X)$ 
is defined on $(x_0, ...., x_{m+n})\in MC_{\e+\delta, m+n}(X)$ by 
\[
\varphi\cdot \psi (x_0, ...., x_{m+n}) = \varphi(x_0, ...., x_{m})\cdot \phi(x_m, ...., x_{m+n})
\]
if $(x_0, ...., x_{m}) \in  MC_{\e, *}(X)$ and $(x_m, ...., x_{m+n}) \in  MC_{\delta, *}(X)$, and by 
\[
\varphi\cdot \psi (x_0, ...., x_{m+n}) = 0
\]
otherwise. The product gives rise to that on the magnitude cohomology $MH_{*}^*(X)$. Moreover, we see that 
the product makes $MH_{*}^*(X)$ into a unital, associative bigraded ring with unit $1=[u]$, where $u \in MC_{0}^0(X)$ is defined by $u(x_0) =1$. Observe that the product is {\it not} graded commutative in general; see \cite[Proposition 2.3]{H}. 

We refer the reader to \cite[Sections 6 and 7]{H} for 
interesting computations of magnitude cohomology rings of finite graphs.

\subsection{The cup product on a cochain complex} 
We consider a commutative diagram 
\begin{eqnarray}\label{eq:F}
\xymatrix@C20pt@R10pt{
F^{\e} \  \ar@{>->}[rr] \ar[rd]_j & & F^{<\e} \ar[ld]^i \ar@{.>}[rdd]_{\alpha : = p\circ i} \ar[r] & F^{<\e}/F^{\e} \ar@{.>}[dd]^{\overline{\alpha}}_\cong\\
 & RP^* \ar[rd]^{p} \ar[ld]_{p'}&  \\
 \K P(<\e)^* & & \K P(\e)^* \ar[ll]^{(\iota_{<\e})^*} & \text{Ker} \ (\iota_{<\e})^* \ar[l]  \ar[r]_{u^*}^\cong & MC_\e^*(X), 
}
\end{eqnarray}
where $F^\e$ and $F^{<\e}$ are kernels of $p$ and $p'$, respectively. 
We have an exact sequence  $0 \to \K P(<\e) \stackrel{\iota_{<\e}}{\longrightarrow} \K P(\e) \to \K P(\e)/\K P(<\e) \to 0$
of chain complexes. 

Let $MC_{\e, *}(X)$ be the unnormalized magnitude chain complex for a generalized metric space $X$; see \cite[Section 4]{H}. 
We observe that there exists an isomorphism 
$u : MC_{\e, *}(X) \to \K P(\e)/\K P(<\e)$ of chain complexes, where $P =N(X) : [0, \infty) \to \mathsf{sSet}$ is the persistence simplicial set described in the beginning of Section \ref{sect:E_nerves}. 
In fact, the inclusion $MC_{\e, *}(X) \to \K P(\e)$ gives rise to the isomorphism, 
Since $\K P(\e)/\K P(<\e)$ is a free $\K$-module,  it follows that the dual sequence 
$$0 \gets \K P(<\e)^* \stackrel{(\iota_{<\e})^*}{\longleftarrow} \K P(\e)^* \gets (\K P(\e)/\K P(<\e))^* \gets 0$$ 
is exact. Then, the natural map $u$ is an isomorphism. 
We have a copersistence cochain complex $F : \V \to \mathsf{coCh_\K}$ defined by $F(\e) = F^{<\e}$ with the natural inclusion 
$F(\e < \e') : F(\e') \to F(\e)$; that is, $F : \V^{\rm op} \to \mathsf{coCh_\K}$ is a functor. 

The usual cup product  $\cup$ on the cochain complex $RP^*:=\text{Hom}_\K(RP, \K)$ is defined by
\begin{eqnarray}\label{eq:cupPD}
(\varphi \cup \psi)(\sigma) = ((d_{p+1}\cdots d_{p+p'})^*\varphi) (\sigma)\cdot_\K ((d_0 \cdots d_{p-1})^*\psi )(\sigma),
\end{eqnarray}
for $\varphi \in (RP)^p$ and $\psi \in (RP)^{p'}$,  where $\cdot_\K$ denotes the multiplication on $\K$. In what follows, we may write $D_+$ and  $D_{-}$ for 
$d_{p+1}\cdots d_{p+p'}$ and $d_0 \cdots d_{p-1}$ in (\ref{eq:cupPD}), respectively.  

Consider the condition (A) below for 
a persistence simplicial set $P: \V \to \mathsf{sSet}$ whose structure maps are injective. 

\begin{itemize}
\item[(A):]  $D_+(\sigma) \in P_m(<\e)$ or  $D_-(\sigma) \in P_n(<\delta)$ for $\sigma \in P_{m+n}(< \e+ \delta)$ and integers $m$, $n$. 
\end{itemize}

\begin{thm}\label{thm:main_I} Let $P$ be a persistence simplicial set for which the condition {\rm (A)} holds. 
 Then, the restriction of the cup product on $RP^*$ to $F(\e)\otimes F(\delta)$ factors through $F(\e+\delta)$. Thus, one has a commutative diagram 
\[
\xymatrix@C25pt@R15pt{
F(\e)\otimes F(\delta) \ar[r]^-\cup \ar[d]_{i\otimes j} & F(\e+\delta) \ar[d]^k\\
RP^* \otimes RP^*  \ar[r]_-\cup  & RP^*,
}
\]
where $i$, $j$ and $k$ are natural inclusions. As a consequence, the product on $F(\ast)$ is associative.  Moreover, the cup product compatible with the structure maps of the copersistence cochain complex; that is, the diagram 
\[
\xymatrix@C25pt@R15pt{
F(\e')\otimes F(\delta') \ar[r]^-\cup \ar[d]_{F(\e\leq \e')\otimes F(\delta\leq \delta')} & F(\e' +\delta') \ar[d]^{F(\e +\delta\leq \e' +\delta')}\\
F(\e)\otimes F(\delta) \ar[r]_-\cup  &  F(\e+\delta)
}
\]
is commutative for $\e\leq \e'$ and $\delta \leq \delta'$. 
\end{thm}

\begin{proof}
By definition, we see that  $F(\e)$ is the kernel of the natural map $RP^* \to \K P(<\e)$. This gives the commutativity of the first square.
The cup product on $RP^*$ is associative and then so is the restriction. 
The direct computation yields the commutativity of the second square. 
\end{proof}

\begin{rem}\label{rem:unit}
Suppose that the unit $1$ in $RP^0$ is in $F(0)^0$. Then $F(\ast)$ is unital in the sense that 
$1\cup \psi = \psi = \psi \cup 1$ for $\psi \in F(\e)$. \todo{Thus, we see that $F$ is a {\it monoidal functor}.}
\end{rem}


\begin{prop}\label{prop:(A)} 
Let $C$ be a $\W^{\rm op}$-category, where $\W$ is the monoidal poset $\V$ or $([0, \infty], \leq, + , 0)$.  
Then, the persistence simplicila set $N(C) : \V \to \mathsf{sSet}$ described in Section \ref{sect:E_nerves} satisfies the condition {\rm (A)}. 
\end{prop}

\begin{proof}
For an element $(x_0, ..., x_m, x_{m+1}, ... , x_{m+n})$ in $N(X)_{n+m}(< \e+\delta)$, we see that $\ell(D_-(x_m, ... ,  x_{n+m})) <\delta$
if $\ell(D_+(x_0, ... , x_{m+n}))\geq \e$. 
\end{proof}

\begin{rem} 
Let $\mathsf{VR}(X) : [0, \infty) \to \mathsf{sSet}$ be the Vietoris-Rips simplicial set of a metric space $(X, d)$ 
defined by 
\[
\mathsf{VR}(X)(\e)_n : = \{(x_0, ..., x_n) \mid x_i \in X, \ d(x_i, x_j) \leq \e \ \text{for} \ i, j  \}.
\] 
In general, the condition (A) for the Vietoris-Rips simplicial set  does {\it not} hold. In fact, for the complete graph $K_n$ with $n\geq 3$, consider the simplex $\sigma =(x_0, x_1, x_2)$ consisting of distinct points. While $D_+(\sigma)=(x_0, x_1)$ and $D_-(\sigma)=(x_1, x_2)$ are not in $P_1(< 1)$, we have $\sigma \in P_2(<2)$. 
\end{rem}

\begin{prop}\label{prop:MC} Let $X$ be a generalized metric space. 
One has a commutative diagram 
\[
\xymatrix@C25pt@R15pt{
F^{<\e}/F^\e \otimes F^{<\delta}/F^\delta \ar[r]^-{\cup} \ar[d]_{\overline{\alpha}\otimes \overline{\alpha}}^\cong& F^{<\e+\delta}/F^{\e+\delta} \ar[d]_{ \overline{\alpha}}^{\cong}\\
\rm{Ker} \ \!(\iota_{<\e})^* \otimes  \rm{Ker} \ \! (\iota_{<\delta})^* \ar[r]^-{\cup} \ar[d]_{u^* \otimes u^*}^\cong
& \rm{Ker} \ \!(\iota_{<\e+\delta})^*   \ar[d]_{u^*}^\cong   \\
 MC_\e^*(X) \otimes  MC_\delta^*(X)  \ar[r]_-{\cdot_H} &  MC_{\e+\delta}^*(X), }
\]
where $\cdot_H$ denotes the product on the (unnormalized) magnitude cochain complex $MC_*^*(X)$
defined by Hepworth; see Section \ref{sect:MH}.  
\end{prop}

\begin{proof}
The result follows from the direct computation.  
\end{proof}

\subsection{The multiplicative magnitude-path spectral sequence (MPSS)}\label{sect:MPSS_1} Let $P : [0, \infty) \to \mathsf{sSet}$ be a persistence simplicial set for which the condition (A) holds. We restrict the domain of the functor $P$ to $\mathbb{N}$. The same notation $P$ is used for the restricted functor. 
We recall  the reachability complex $RP:=\text{colim}_{l\in \mathbb{N}}\K P(l)\cong \K(\text{colim}_{l\in \mathbb{N}}P)$. 

With the same notation as in the diagram (\ref{eq:F}), we have the decreasing filtration 
$F=\{F^l(RP^*)\}_{l \in \mathbb{N}}$ defined by 
\[
(F^l(RP^*))^n := (F^{<l})^n = (F^{l-1})^n.
\]
The usual procedure to make a spectral sequence, we have a right half-plane spectral sequence $\{E_r^{*,*}, d_r\}$. 
In fact, for a filtered complex $(A, F, d)$, 
the $E_r$-term of the spectral sequence $\{E_r^{*, *}, d_r\}$ is defined  by 
$E_r^{p, q} := Z_r^{p, q} / (Z_{r-1}^{p+1, q-1}+ B^{p, q}_{r-1})$, 
where 
$Z_r^{p,q}$ and
$B_r^{p,q}$ defined by
\begin{equation}\label{eq:Z-B}
Z_r^{p, q} := F^pA^{p+q}\cap d^{-1}(F^{p-r}A^{p+q+1})  \  \ \text{and} \ \
B_r^{p, q}:= F^pA^{p+q}\cap d(F^{p-r}A^{p+q-1}), 
\end{equation}
respectively. 
The differential $d_r$ on the $E_r$-term is induced by $d$; see \cite[The proof of Theorem 2.6]{M}. 
Thanks to Theorem \ref{thm:main_I}, the spectral sequence is endowed with a product structure. 

The following theorem is a corollary of the consideration above. 
\begin{thm}\label{thm:SS} \rm{(The multiplicative magnitude-path spectral sequence)} 
Let $X$ be a generalized metric space. Then, the filtration $F:=\{F^l\K N(X)\}_{l\in \mathbb{N}}$ gives rise to a multiplicative spectral sequence 
$\{E_r^{*,*}, d_r\}_{r\geq 1}$ with $E_r^{p, q} = 0$ for $p+q < 0$,  $E_1^{p, q}\cong MH_p^{p+q}(X)$ as a $\K$-module and 
\[
E_1^{*, *}\cong \oplus_{p, q}MH_p^{p+q}(X) 
\] 
as a bigraded algebra. Moreover, the spectral sequence converges to $H^*(RN(X))$ as an algebra provided the filtration $F$ is bounded.
\end{thm}

\begin{proof} 
It follows from  Proposition \ref{prop:(A)} and Theorem \ref{thm:main_I} that \[
\cup( F^sN(X)^m\otimes F^t N(X)^n)  \subset F^{s+t}N(X)^{m+n}
\]
for the cup product.  
Then, we see that the filtration $F$ gives a bigraded algebra. 
By \cite[Theorem 2.14]{M}, we have the multiplicative spectral sequence associated by $F$. 
Thus, Proposition \ref{prop:MC} enables us to deduce that the $E_1$-term is isomorphic to the magnitude cohomology ring 
$\oplus_{l, n}MH_l^n(X)$.  We have the result. 
 \end{proof}
 
 \begin{rem} The $p$-axis $E_1^{*,0}$ in the spectral sequence in Theorem \ref{thm:SS} is nothing but 
the diagonal part of the magnitude cohomology ring $\oplus_{n, n}MH_n^{n}(X)$. 
 \end{rem}
  
 \begin{rem} i) If the generalized metric space $X$ is finite, then the filtration $F$  in Theorem \ref{thm:SS} is bounded.
 
ii) Suppose that the filtration $F$ is {\it degreewise bounded}; that is, for each cohomology degree $n$, there exists a filtration degree $l=l(n)$ such that $F^l\K P^n = RP^n$. 
 Then, for $(p, q)$ a pair of nonnegative integers,  there exists an integer $N$ such that $E_r^{p, q} \cong E_\infty^{p, q}$ as a $\K$-module for $r \geq N$;  see \cite[The proof of Theorem 2.6]{M}. 
 
iii) We may apply the construction of the spectral sequence in the proof of Theorem \ref{thm:SS} to the enriched nerve of the $\overline{\mathbb{N}}^{\rm op}$-category $\lceil \  \rceil X$; see Example \ref{ex:L} ii). The $E_1$-term is isomorphic to the magnitude cohomology algebra of $\lceil \  \rceil X$. 

iv) Let $X$ be a compact Riemannian manifold, more general, a Riemannian diffeological space which is compact for the D-topology; see \cite[Theorem 3.21]{KSS}.
Then the map $\ell : X^n \to [0, \infty]$ is bounded for each $n$; see \cite[Theorem 3.21]{KSS} for a diffeological space. This implies that  the filtration $F$ which construct the spectral sequence 
in Theorem \ref{thm:SS} and ii) above is bounded degreewise. Then, these spectral sequences converge degreewise. 
 \end{rem}
 
\begin{rem} 
The cohomology ring $H^*(RN(X))$ is commutative. In fact, the cochain complex $RN(X)$ is obtained 
by the simplicila set $\text{colim}_{l\in \mathbb{N}}N(X)$. 
However, the magnitude cohomology ring is not necessarily commutative; see \cite[Proposition 2.3]{H}. 
It is an interesting problem to consider which term in the multiplicative MPSS is (bigraded) commutative. 
\end{rem}

We end this section with a comment on a relationship between the MPSS and an $r$-homotopy equivalence mentioned in Section \ref{sect:E_nerves}. The topic is already discussed in \cite[Section 4]{Asao23}, \cite[Section 3]{LRW} and \cite[Section 2]{HR_2026}. 

Let $X$ and $X'$ be generalized metric spaces and $f, g :  X \to X'$ $r$-homotopic maps. Assume that $r$ is a non-negative integer. 
The simplicial homotopy $h : N(X) \to N(X')$ from $N(f)$ to $N(g)$ in Remark \ref{rem:homotopy} gives rise to a chain homotopy 
$\widetilde{h} : RN(X)_* \to RN(X')_{*+1}$ from $RN(f)$ to $RN(g)$; see, for example, \cite[Chapter 3, 4.2 Proposition]{Barr}. 
By the construction,  we see that  $\widetilde{h}(\K N(X)(p))\subset \K N(X')(p+r)$ for each $p \in N$. This enables us to conclude that $\widetilde{h}^* (F^{p+r}(RN(X')^*)) \subset F^{p}(RN(X)^{*-1})$. Thus, by \cite[Proposition XV, 3.1]{C-E}, 
we have 

\begin{prop}\label{prop:r-H_MPSS} {\rm (\cite[Theorem 4.37]{Asao23} and \cite[Proposition 3.3.3]{LRW})}
Suppose that $f, g : X \to X'$ are $r$-homotopic maps of generalized metric spaces. Then, $E_s(f) = E_s(g)$ for all $s > r$. As a consequence, for an $r$-homotopy equivalence $f : X \to X'$, the morphism $E_r(f) : E_r^{*,*}(X) \to E_r^{*,*}(X')$ is a quasi-isomorphism of differential bigraded algebras.  
\end{prop}

\section{Interleavings and the MPSS}\label{sect:ID_MPSS}

We describe  a relationship between interleavings of enriched nerves and the notion of an $r$-homotopy 
in Theorem \ref{thm:interleavings}. Proposition \ref{prop:r-H_MPSS} relates the $r$-homotopy to the MPSS. 
Therefore, one may expect a relationship of the MPSS and the interleaving between persistence chain complexes. 
However, in general, an interleaving is not a morphism of persistence objects and, therefore, does not naturally induce  
a morphism of spectral sequences. 
The following proposition gives a relationships of the MPSS and interleaving with an acyclic persistence chain complex.

\begin{thm}\label{thm:IHC_MPSS}
Let $X$ be a generalized metric space whose distance takes values in non-negative integers and $\{E_r^{*,*}(X), d_r\}$ the MPSS for $X$. Suppose that 
\[
\displaystyle d_{\rm IHC}^{\rm \ \! et}(\K N(X), \K N(pt)) = \frac{r}{2},
\] 
where $d_{\rm IHC}^{\rm \ \!et}$ denotes the interleaving distance in the homotopy category of the category $\mathsf{et}(\mathsf{Ch}_{\K}^{[0, \infty)})$ of extended tame persistence chain complexes endowed with the model structure described in \cite{KNSWY}; see Theorem \ref{thm:modelST} below for the detail. 
Then, one has  $E_s^{*,*}(X) =\K$ for $s >r$; see also Remark \ref{rem:homology}. 
In particular, the reachability complex 
$RN(X)= {\rm colim}_l \K N(X)(l)$ is acyclic. 
\end{thm}

Before proving Theorem \ref{thm:IHC_MPSS}, we need to explain terminology used in the theorem. 
The model structure mentioned in Theorem \ref{thm:IHC_MPSS} is given modifying that of the category of tame persistence objects defined in \cite{CGL}.  The advantage of using the model structure is that we have comparatively tractable cofibrations in the model category. On the other hand, it is difficult to say that a cofibration in the projective model structure of a functor category has an explicit form; see, for example, \cite[Theorem 11.6.1]{H_Model}.

We here consider the persistence chain complex $\K N(X)$ in the category $\mathsf{et}(\mathsf{Ch}_{\K}^{[0, \infty)})$ of extended tame persistence chain complexes.

\begin{defn}\label{defn:et} {\rm (\cite[Definition 2.1]{CGL}, \cite[Definition 2.4]{KNSWY})}
Let $\C$ be a category. 
A sequence $\tau_0 < \tau_1 < \cdots  <\tau_n<\cdots$ in $[0,\infty)$, which is divergent or finite,
{\it discretises} a functor $X: [0, \infty) \to \C$ if $X({s\leq t}):X(s)\to X(t)$ may fail to be an isomorphism only if there is $a\in \mathbb{N}$ such that $s<\tau_a\leq t$. A functor $X : [0, \infty)\to \C$ is called {\it extended tame} if there is a sequence that discretises it. Let  $\et{\C}$ denote the full subcategory of the functor category ${\C}^{[0, \infty)}$ whose objects are extended tame functors. 
An object $X$ in $\et{\C}$ is called {\it tame} if the sequence which discretises $X$ is finite; see \cite{CGL}.
\end{defn}

Let $\C$ be a category. 
We recall a factorization of a morphism $g : X\to Y$ in $\et{\C}$. Let  $\tau_0 < \tau_1 < \cdots  <\tau_n<\cdots $ be a sequence discretising both $X$ and $Y$; see \cite[Section 2.2]{CGL}.  By induction on
$\{0, 1, 2, ..., n, ... \}$, we define morphisms
$\bar{g}({\tau_a}) : X({\tau_a}) \to Q({\tau_a})$ and $\hat{g}({\tau_a}) : Q({\tau_a}) \to Y({\tau_a})$ in $\C$ as follows:
%
$(\bar{g}(0) : X(0) \to Q(0)):= (1 : X(0) \to X(0))$ \ and \  $(\widehat{g}(0) : Q(0) \to Y(0)):=(g(0) : X(0)\to Y(0))$.
For $a > 0$, the object $Q({\tau_a})$ is defined by 
\[\text{colim}(\xymatrix@C35pt@R20pt{Y({\tau_{a-1}}) & X({\tau_{a-1}}) \ar[l]_-{g({\tau_{a-1}})} \ar[r]^-{X({\tau_{a-1} < \tau_a})} & X({\tau_a}) })
\]
with $\bar{g}({\tau_a}) : X({\tau_a}) \to Q({\tau_a})$ and $\widehat{g}({\tau_a}) : Q({\tau_a}) \to Y({\tau_a})$ which fit in the commutative diagram
\begin{eqnarray}\label{diagram:the_pushout}
\xymatrix@C40pt@R15pt{
 X({\tau_{a-1}}) \ar[r]^{X({\tau_{a-1} < \tau_a})} \ar[d]_{g({\tau_{a-1}})}& X({\tau_a}) \ar[d]^{\bar{g}({\tau_a})}  \ar@/^0.8pc/[rdd]^{g({\tau_a})} &\\
 Y({\tau_{a-1}}) \ar@/_0.8pc/[rrd]_{Y({\tau_{a-1} < \tau_a})}\ar[r] & Q({\tau_a}) \ar[rd]^{\widehat{g}({\tau_a})} & \\ 
  & & Y({\tau_a})
}
\end{eqnarray}
consisting of the pushout square.
For $a > 0$, define $Q({\tau_{a-1}< \tau_a}) : Q({\tau_{a-1}}) \to Q({\tau_a})$ to be the composite of the morphisms $\widehat{g}({\tau_{a-1}}) : Q({\tau_{a-1}})  \to Y({\tau_{a-1}})$ and $Y({\tau_{a-1}}) \to Q({\tau_a})$ which is the bottom horizontal arrow in (\ref{diagram:the_pushout}). 
Thus, we have a factorization $g = \widehat{g}\circ \bar{g}$ of $g$.

Let $\mathcal{M}$ be a model category. 
We recall a model structure on the full subcategory of the funcotr category  $\mathcal{M}^{[0, \infty)}$ consisting of extended tame persistence objects.   
The following results for {\it tame} objects are originally proved by  Chach\'olski, Giunti and Landi \cite{CGL}. The extended tame versions are  also proved by the same argument as in the proof of \cite[Theorem 2.2]{CGL}; see \cite[Comments after Proposition 2.6]{KNSWY}.

\begin{thm}\label{thm:modelST}{\rm (\cite[Theorem 2.2]{CGL}, \cite[Theorem 2.5]{KNSWY})} Let $\modelcat$ be a model category.  
With the notation and terminology mentioned above, the following choices of weak equivalences, fibrations and cofibrations form a model structure on $\et{\modelcat}$. A morphism $g:X\to Y$ in $\et{\modelcat}$ is a
    \begin{itemize}
        \item weak equivalence if $g(t):X(t)\to Y(t)$ is a weak equivalence for all t,
        \item fibration if $g(t):X(t)\to Y(t)$ is a fibration for all t,
        \item cofibration if $\widehat{g}(t): Q(t)\to Y(t)$ is a cofibration for all t.
    \end{itemize}
\end{thm}

It is worthwhile mentioning that the model category in Theorem  \ref{thm:modelST} is not necessarily cofibrantly generated. 
We have a characterization of a cofibrant object in $\et{\modelcat}$.

\begin{prop}\label{prop:cofibrantOb}{\rm (\cite[Proposition 2.3]{CGL}, \cite[Proposition 2.6]{KNSWY})}
An object $X$ in $\et{\modelcat}$ is cofibrant if and only if $X(0)$ is cofibrant and, for any $s < t$ in $[0, \infty)$, the transition morphism $X({s<t}) : X(s) \to X(t)$ is a cofibration in $\modelcat$.
\end{prop}

We apply the construction in Theorem \ref{thm:modelST} to the category $\mathsf{Ch_\K}$ of non-negatively graded chain complexes whose model structure is described in \cite[7.2 Theorem]{DS}. Observe that by definition, a morphism $i : M \to N$ in $\mathsf{Ch_\K}$ is a cofibration if $i_k : M_k \to N_k$ is a monomorphism and the cokernel of $i_k$ is a projective $\K$-module for $k \geq 0$. 
 
Let $I$ be the chain complex over $\K$ generated by $\widetilde{0}$, $\widetilde{1}$ and $\iota$ with $d(\iota) = \widetilde{0} -\widetilde{1}$ and $\deg(\iota) = 1$. We regard $I$ as the constant persistence chain complex.

Let $M$ be a cofibrant object in $\et{\mathsf{Ch}_\K}$. For the folding map $\nabla: M^{\oplus 2} \to M$, 
we have a commutative diagram
 \begin{eqnarray}\label{diagram:factrisation_I}
\xymatrix@C40pt@R12pt{
M^{\oplus 2} \ar@/^1.5pc/[rr]^{\nabla} \ar[r]_-{g} & M\otimes I \ar@{->>}[r]_-\pi^-\sim & M
}
\end{eqnarray}
in $\et{\mathsf{Ch}_\K}$. Here $M\otimes I$ is the persistence chain complex defined by the tensor product pointwise. 
Moreover, the map $g$ is  induced by the natural inclusion and the map $\pi$ is  defined by $\pi(x\otimes \widetilde{i}) = x$ for $i = 0, 1$ and $\pi(x\otimes \iota)=0$. It is readily seen that $\pi$ is a trivial fibration. 

\begin{lem}\label{lem:a_cylinderOb} Let $M$ be a cofibrant object in $\et{\mathsf{Ch}_\K}$. The map $g$ in {\rm (\ref{diagram:factrisation_I})} is a cofibration in the sense in Theorem \ref{thm:modelST}.  As a consequence, 
the persistence chain complex $M\otimes I$ is a cylinder object of $M$ in $\mathsf{et}(\mathsf{Ch}_{\K}^{[0, \infty)})$. 
\end{lem}

\begin{proof} It sufficeis to show that $g$ in the diagram (\ref{diagram:factrisation_I}) is a cofibration in $\mathsf{et}(\mathsf{Ch}_{\K}^{[0, \infty)})$. 
We consider the commutative  diagram 
\begin{eqnarray*}
\xymatrix@C30pt@R15pt{
M(k)^{\oplus 2}  \ar[r]^{\eta^{\oplus 2}({k <  k+1})} \ar[d]_{g({k})}& 
M(k+1)^{\oplus 2} \ar[d]^{\bar{g}({k+1})}  \ar@/^7.0pc/[dd]^{g({k+1})} &\\
M(k) \otimes I  \ar@/_0.8pc/[rd]_-{(\eta\otimes I)({k<k+1})\ \ \ } \ar[r] &
M(k+1)^{\oplus 2} \oplus  M(k)\otimes \K\{\iota\} \ar[d]^{\widehat{g}(k+1)}& \\ 
  & M(k+1)\otimes I &
}
\end{eqnarray*}
in which the inner square is the pushout.  It follows that $\widehat{g}(k+1)$ is a monomorphism in each degree of the chain complex and the cokernel of $\widehat{g}(k+1)$ is a projective module in each degree. Then, the morphism 
$\widehat{g}(k+1)$ is a cofibration in $\mathsf{Ch}_{\K}$; see the comment that follow Proposition \ref{prop:cofibrantOb}. Thus, by definition, we see that $g$ is a cofibration in $\mathsf{et}(\mathsf{Ch}_{\K}^{[0, \infty)})$.
\end{proof}

\begin{proof}[Proof of Theorem \ref{thm:IHC_MPSS}]
Let $\K$ be the constant persistence chain complex $\K$ with $\K(l)_0=\K$ and $\K(l)_i = 0$ for $i\geq 1$. 
Since $\K N(pt) \cong \K$ in the homotopy category ${\rm Ho}(\mathsf{et}(\mathsf{Ch}_{\K}^{[0, \infty)}))$, it follows form the assumption that 
$d_{\rm IHC}^{\rm \ \! et}(\K N(X), \K) = r/2$. 

We observe that the non-negative integers discretises the persistence chain complex $\K N(X)$. Therefore, we see that $r$ is an integer. In fact, if $r$ is not an integer, then $\K N(X)$ and $\K$ are 
$(\lfloor r \rfloor/2)$-interleaved in the homotopy category, which is a contradiction. 

Proposition \ref{prop:cofibrantOb} yields that $\K N(X)$ and $\K$ are cofibrant.  Moreover, the both of the objects are fibrant.  Thus, we have an $r/2$-interleaving $(\K N(X), \K, \varphi, \psi)$ in the homotopy category for which 
$\psi^{r/2}\circ \varphi \simeq \eta^{r}_{N(X)}$ and $\varphi^{r/2}\circ \psi \simeq \eta^{r}_\K = id_K$ via homotopies 
$H_1 :  \K N(X)\otimes I \to \K N(X)^r$ and $H_2 : \K\otimes I \to \K^r$ with the cylinder objects described in 
Lemma \ref{lem:a_cylinderOb}, respectively; see Definition \ref{defn:interleaving}.

Since $\K$ is concentrated in the zero degree of each chain complex, it follows that 
$\varphi : \K N(X) \to \K^{r/2}$ factors through $\K$; that is, $\varphi = \eta^{r/2}_\K\circ \varphi'$, for some morphism 
$\varphi' : \K N(X) \to \K$ of persistence chain complexes.  Moreover, the image of $\psi$ is in $\K N(X)(\ast)_0$ and then 
in the image of the structure map $\eta : \K N(X)(0)_0 \to\K N(X)(\ast)_0$. Thus, there exists a morphism 
$\psi' : \K  \to \K N(X)$ of persistence chain complexes such that $\psi = \eta^{r/2}_{N(X)}\circ \psi'$.

The homotopy $H_1$ mentioned above is in $\mathsf{et}(\mathsf{Ch}_{\K}^{[0, \infty)})$. Therefore, the homotopy induces a chain map $RN(X)\otimes I \to RN(X)$ and then a chain homotopy $h_1 : RN(X)_* \to RN(X)_{*+1}$ defined by 
$h_1(x) = (-1)^{\deg x}H_1(x\otimes \iota)$. We have $h_1 : \psi'\varphi' \simeq id_{RN(X)}$. The same argument above yields 
a chain homotopy $h_2$ from $\K = {\rm colim} \K$ to itself which gives  $\varphi' \psi'\simeq id_\K$. The homotopies $h_i$ increases the filtration degree of $RN(X)$ and $\K$ by $s$. Thus, the argument before Proposition \ref{prop:r-H_MPSS} yields that $E_s^{*,*} = \K$ for $s > r$. 
\end{proof}

\begin{rem}\label{rem:homology} Let $X$ be a generalized metric space which satisfies the same condition as in Theorem \ref{thm:IHC_MPSS}. Consider the homological MPSS $\{E^r_{*,*}(X), d^r\}$ for $X$. Then, 
the argument before Proposition \ref{prop:r-H_MPSS} and the proof of Theorem \ref{thm:IHC_MPSS} 
yield that $E^s_{*, *}=\K$ for $s >r$. 
\end{rem}

\begin{rem} i) 
Let $G$ be a connected finite digraph with the shortest path metric. Then $cN(G)$ is contractible and then $RN(G)$ is acyclic; 
see Remark \ref{rem:homotopy} (2). Therefore, the distance  $d_{\rm IHC}^{\rm \ \! et}(\K N(X), \K N(pt))$ is finite. 

ii) The distance in the homotopy category described in Theorem \ref{thm:IHC_MPSS} is computed by using the bottleneck distance between the persistence homology $H_*(\K N(X))$ and the constant persistence module $\K$ if $\K$ is a field; see Remark \ref{rem:et_CohID} for more details. 
\end{rem}

 The interleaving distance $d_{\rm IHC}^{\rm \ \!et}$ described in Theorem \ref{thm:IHC_MPSS} is related to $d_{\rm IHC}$ considered in Section \ref{sect:E_nerves}; see Appendix \ref{sect:et_proj}. The relationship enables us to 
 give a computational example of the interleaving distance in the homotopy category between enriched nerves 
of a direct cycle and the point. 
We recall the definition of the directed cycle. 

\begin{defn}{\rm (\cite[Definition 8.1]{HR})} For $m\geq 1$, the {\it directed cycle} $Z_m$ of length $m$ is the digraph with $m$ cyclicially ordered vertices, and with a single directed edge adjacent vertices, consistently oriented. 
\end{defn}

In \cite[Section 8]{HR}, the homological MPSS is explored in detail 
the directed cycle $Z_m$. By using the result of the computation, we have 

\begin{prop}\label{ex:DirectedCycles} For the directed cycle $Z_m$ of length $m \geq 3$, 
\[
\displaystyle d_{\rm IHC}^{\rm \ \! et}(\K N(Z_m), \K N(pt)) = d_{\rm IHC}(N(Z_m)_\bullet, N(x_1)_\bullet) = \frac{m-1}{2}.
\]
\end{prop}

\begin{proof}
By virtue of Corollary \ref{cor:diameter}, we have $d_{\rm IHC}(N(Z_m)_\bullet, N(x_1)_\bullet) \leq  \frac{m-1}{2}$. It follows from 
Corollary \ref{cor:N(X)} that $d_{\rm IHC}^{\rm \ \! et}(\K N(Z_m), \K N(pt)) \leq d_{\rm IHC}(N(Z_m)_\bullet, N(x_1)_\bullet)$.  
The result \cite[Theorem 8.2]{HR} enables us to deduce that for $\{E^r_{*,*}(Z_m), d^r\}$ the homology MPSS of $Z_m$, the $E_m$-term is trivial but $E_{m-1}$-term is non-trivial. Then, Theorem \ref{thm:IHC_MPSS} and Remark \ref{rem:homology}
imply that $\frac{m-1}{2} \leq d_{\rm IHC}^{\rm \ \! et}(\K N(Z_m), \K N(pt))$. We have the result.  
\end{proof}

\medskip
\noindent
{\it Acknowledgements.}
The author thanks Takahito Naito, Shun Wakatsuki and Toshihiro Yamaguchi for valuable discussions on interleavings of enriched nerves. The author was partially supported by JSPS KAKENHI
Grant Number JP26K00604.

\appendix

\section{$d_{\rm IHC}$ versus $d_{\rm IHC}^{\rm \  \! et}$}\label{sect:et_proj}
We compare the interleaving distances $d_{\rm IHC}$ and $d_{\rm IHC}^{\rm \  \! et}$ between enriched nerves. Let $\mathcal{M}$ be a cofibrantly generated model category and $i : \mathsf{et}(\mathcal{M}^{[0, \infty)}) \to \mathcal{M}^{[0, \infty)}$ the inclusion functor. 
We regard $\mathcal{M}^{[0, \infty)}$ and  $\mathsf{et}(\mathcal{M}^{[0, \infty)})$ as 
the model categories endowed with the projective model structure and the model structure described in Theorem \ref{thm:modelST}, respectively. 

Let $(\lfloor  \  \rfloor)_* : \mathcal{M}^{{\mathbb{N}\cup \{0\}}} \to \mathsf{et}(\mathcal{M}^{[0, \infty)})$ be the functor induced by the floor function $\lfloor  \  \rfloor : [0, \infty) \to {\mathbb{N}\cup \{0\}}$. Observe that the image of the functor $(\lfloor  \  \rfloor)_*$ is in $\mathsf{et}(\mathcal{M}^{[0, \infty)})$. Thus, we have $(\lfloor  \  \rfloor)_* = i \circ (\lfloor  \  \rfloor)_*^|$, where $(\lfloor  \  \rfloor)_*^| :   \mathcal{M}^{{\mathbb{N}\cup \{0\}}} \to \mathsf{et}(\mathcal{M}^{[0, \infty)})$ is the corestriction. 

\begin{prop}\label{prop:et_proj}
For objects $F$ and $G$ in the image of the functor $(\lfloor  \  \rfloor)_*^|$ , one has 
\[
d_{\rm IHC}^{\rm \ \! et}(F, G) = d_{\rm IHC}(iF, iG),
\]
where $d_{\rm IHC}^{\rm \ \! et}$ denotes the interleaving distance in the homotopy category of $\mathsf{et}(\mathcal{M}^{[0, \infty)})$. 
\end{prop}

\begin{proof}
Let $res :  \mathcal{M}^{[0, \infty)} \to \mathsf{et}(\mathcal{M}^{[0, \infty)})$ be the composite of $(\lfloor  \  \rfloor)_*^|$ and  functor induced by the inclusion $\iota : \mathbb{N}\cup \{0\} \to [0, \infty)$ of posets. 
We observe that $res \  i F = F$ for the object in the image of the functor $(\lfloor  \  \rfloor)_*^|$. 
The functors $i$ and $res$ preserve weak equivalences. Then, we have the right derived functors ${\bf R}i :  {\rm Ho}(\mathsf{et}(\mathcal{M}^{[0, \infty)})) \to {\rm Ho}(\mathcal{M}^{[0, \infty)})$ and ${\bf R}res :   {\rm Ho}(\mathcal{M}^{[0, \infty)}) \to {\rm Ho}(\mathsf{et}(\mathcal{M}^{[0, \infty)}))$ between homotopy categories.  
Since each object in $\mathcal{M}^{[0, \infty)}$ and  $\mathsf{et}(\mathcal{M}^{[0, \infty)})$ is fibrant, it follows from \cite[Proposition 9.3]{DS} that 
${\bf R} iF \cong iF$ and ${\bf R} res \ \!G \cong G$ for $F \in  \mathsf{et}(\mathcal{M}^{[0, \infty)})$ and $G \in \mathcal{M}^{[0, \infty)}$. 
Thus, we have 
\begin{align}
d_{\rm IHC}^{\rm \ \! et}(F, G) &\geq d_{\rm IHC}({\bf R} iF, {\bf R} iG)= d_{\rm IHC}(iF, iG) \\
& \geq d_{\rm IHC}^{\rm \ \! et}({\bf R}res\  iF, {\bf R} res \ iG) =  d_{\rm IHC}^{\rm \ \! et}(res \ iF, res \ iG)  \\
&=  d_{\rm IHC}^{\rm \ \! et}(F, G). 
\end{align}
This completes the proof
\end{proof}

\begin{cor}\label{cor:N(X)} For generalized metric spaces $X$ and $Y$ whose distances take values in non-negative integers, one has 
\[
d_{\rm IHC}^{\rm \ \! et}(N(X), N(Y)) = d_{\rm IHC}(N(X), N(Y)), 
\]
\[
d_{\rm IHC}^{\rm \ \! et}(\K N(X), \K N(Y)) = d_{\rm IHC}(\K N(X), \K N(Y)), 
\]
\[
d_{\rm IHC}^{\rm \ \! et}(\K N(X), \K N(Y)) \leq d_{\rm IHC}^{\rm \ \! et}(N(X), N(Y)) \  {\rm and} 
\]
\[
d_{\rm IHC}(\K N(X), \K N(Y)) \leq d_{\rm IHC}(N(X), N(Y)).
\] 
\end{cor}

\begin{proof} With the same notation as in the proof of Proposition \ref{prop:et_proj}, 
we see that $N(X) = res\ i N(X)$ in $\mathsf{et}(\mathsf{sSet}^{[0, \infty)})$ and 
 $\K N(X) = res\ i \K N(X)$ in $\mathsf{et}(\mathsf{Ch}_\K)^{[0, \infty)}$. The equalities follow from Proposition \ref{prop:et_proj}. 
By virtue of \cite[Proposition 3.6]{B-S}, we have the inequalities.
\end{proof}

We further relate the distance $d_{\rm IHC}(\K N(X), \K N(Y))$ to the {\it cohomology interleaving distance} $d_{\rm CohI}$ defined in \cite[Definition 3.1]{KNWY}. Write $d_{\rm IHC}(F, G)_{\mathcal{C}}$  for the interleaving distance in the homotopy category between objects $F$ and $G$ in the category $\mathcal{C}= \mathsf{Ch_\K}^{[0, \infty)}$ or $\mathsf{Ch_\K}^\R$. 
Moreover, we denote by $d_{\rm I}(H, K)_{\mathcal{D}}$  the interleaving distance between objects $H$ and $K$ in $\mathcal{D} = \mathsf{Mod_\K}^{[0, \infty)}$ or 
$\mathcal{D} = \mathsf{Mod_\K}^\R$ the category of persistence modules over $\K$. 
In particular, for $X$ and $Y$ in Corollary \ref{cor:N(X)},  the distance $d_{\rm IHC}(\K N(X), \K N(X))_{\mathsf{Ch}_\K^{[0, \infty)}}$ 
is nothing but 
$d_{\rm IHC}(\K N(X), \K N(Y))$ 
 with the notation in Section \ref{sect:I-r} and Corollary \ref{cor:N(X)}. 

\begin{prop}\label{prop:ChoID} Let $\K$ be a field. For objects $F$ and $G$ in $\mathsf{Ch}_\K^{[0, \infty)}$, one has 
\[
d_{\rm IHC}(F, G)_{\mathsf{Ch}_\K^{[0, \infty)}} = d_{\rm ChoI}(F, G) := {\rm sup}\big\{d_{\rm I}(H_n(F), H_n(G))_{\mathsf{Mod_\K^{[0, \infty)}}} \mid n \in {\mathbb Z} \big\}, 
\]
where $H_n$ is the $n$th homology functor form $\mathsf{Ch}_\K$  to $\mathsf{Mod}_\K$ the category of $\K$-modules. 
\end{prop}

\begin{proof}
Consider the maps 
$\xymatrix@C18pt@R20pt{
[0, \infty) \ar[r]^-j & \R \ar[r]^-k & [0, \infty)
}
$
defined by the inclusion and $k(x) = 0$ for $x\leq 0$ and $k(x) = x$ for $x\geq 0$. Then, the maps induce functors 
\[
\xymatrix@C18pt@R20pt{
\mathcal{D}^{[0, \infty)} & \mathcal{D}^\R  \ar[l]_-{j^*}& \mathcal{D}^{[0, \infty)} \ar[l]_-{k^*} 
}
\]
with $j^*k^*=id$ for $\mathcal{D} = \mathsf{Ch_\K}$ and $\mathsf{Mod}_\K$. 
We have a sequence of equalities 
\begin{align}
d_{\rm IHC}(F, G)_{\mathsf{Ch}_\K^{[0, \infty)}}  &= d_{\rm IHC}(k^*F, k^*G)_{\mathsf{Ch}_\K^\R}  \\
& \ \ \ \  \text{(by the same argument as in the proof of Proposition \ref{prop:et_proj}} \\
  &  \ \ \ \ \  \text{with the right derived functors of  $k^*$ and $j^*$)}  \\
  & = {\rm sup}\big\{d_{\rm I}(H_n(k^*F), H_n(k^*G))_{\mathsf{Mod_\K^{\R}}} \mid n \in {\mathbb Z} \big\} \ 
  \text{(by \cite[Theorem 3.3]{KNWY})} \\
  & = {\rm sup}\big\{d_{\rm I}(k^*H_n(F), k^*H_n(G))_{\mathsf{Mod_\K^{\R}}} \mid n \in {\mathbb Z} \big\} \\ 
  & \ \ \ \  
  \text{(by the commutativity of $k^*$ and $H_l$)} \\
  &= {\rm sup}\big\{d_{\rm I}(H_n(F), H_n(G))_{\mathsf{Mod_\K^{[0, \infty)}}}\mid n \in {\mathbb Z} \big\} 
  \   \text{(by $j^*k^* = id$)} \\
  &=  d_{\rm ChoI}(F, G). 
\end{align}
We need the assumption that $\K$ is a field to apply \cite[Theorem 3.3]{KNWY} in the argument above. 
The fact that the polynomial algebra $\K[t]$ is a {\it graded} principal ideal domain is the key to proving the theorem; see 
\cite[Lemmas 3.9 and 3.10]{KNWY} for the details.  
\end{proof}

\begin{rem}\label{rem:et_CohID} Suppose that the underlying ring $K$ is a field. Then, 
we have 
\[
d_{\rm IHC}^{\rm \ \! et}(\K N(X), \K N(Y)) = d_{\rm IHC}(\K N(X), \K N(Y)) = d_{\rm ChoI}(\K N(X), \K N(Y))
\] 
for generalized metric spaces $X$ and $Y$ whose distances take values in non-negative integers. 
This follows from Corollary \ref{cor:N(X)}  and Proposition \ref{prop:ChoID}. 

The interleaving distance of persistence $\K$-module $H$ and $K$, which are of finite dimension pointwise, 
is determined by the bottleneck distance between $H$ and $K$ via the interval module decompositions of the modules; see, for example,  \cite[Theorem 4.16]{B-S}.  
Thus, in particular, the interleaving distance $d_{\rm IHC}^{\rm \ \! et}(\K N(X), \K N(pt))$ in 
Theorem \ref{thm:IHC_MPSS} may be computed with the bottleneck distance between $H_n(N(X); \K)$ and the constant persistence module $\K$ for each $n$. 
\end{rem}

\end{document}